%

\input amstex
 
\documentstyle{amsppt}
\topmatter
\title A Note on UMD Spaces and Transference in
 Vector-valued Function Spaces\endtitle
 
\author Nakhl\'e H. Asmar,  Brian P. Kelly, and Stephen
Montgomery-Smith \endauthor
 
\thanks The work of the first and third authors was
partially funded by NSF grants.
The second and third authors' work was partially funded by the
University of Missouri Research Board\endthanks

 
\abstract  A Banach space $X$ is called an HT 
space if the Hilbert transform is bounded from
$L^p(X)$ into $L^p(X)$, where $1<p<\infty$.
We introduce the notion of an ACF Banach space, 
that is, a
Banach space $X$ for which we have
an abstract M. Riesz Theorem for
conjugate functions in $L^p(X)$, $1<p<\infty$.
Berkson, Gillespie, and Muhly
\cite{5} showed that
$X\in\text{HT}\implies
X\in\text{ACF}$.
In this note, we will show
that $X\in\text{ACF}\implies
X\in\text{UMD}$, thus providing a new
proof of Bourgain's result 
$X\in\text{HT}\implies
X\in\text{UMD}$.
\endabstract
 
\rightheadtext{UMD Spaces and Vector-valued Transference}
\endtopmatter
 
\document
 
\def\mysum{\sum\limits}

\def\myint{\int\limits}
 
\def\msp{\Omega}
\def\malg{ {\Cal F} }

\def\lqq{\lq\lq}
\def\rqq{\rq\rq}
 
\def\C{\Bbb C}
\def\N{\Bbb N}
\def\R{\Bbb R}
\def\T{\Bbb T}
\def\Z{\Bbb Z}
\def\dlgrp{\Gamma}
\def\dualz{{\Z^{\infty\,\ast}}}

\def\od{P}

\def\UMD{\text{UMD}}

\def\sgn{\text{sgn}}

\def\LBspc#1{L^{#1}(X)}

\def\wprod{ {\prod\limits_{n\in \N}}\hskip-
2pt\vphantom{{\prod\limits_{n\in \N}}}^\ast\,\Z}

\def\endprf{\hfill $\blacksquare$}
 
\loadmsbm
\loadeusm
\magnification=\magstep1
\baselineskip=21pt
\CenteredTagsOnSplits

\specialhead 1. Introduction \endspecialhead
 
Recent structure theorems of orders and results in abstract
harmonic analysis on groups with ordered dual groups have
shown that abstract harmonic analysis on such groups
captures martingale theory and classical 
harmonic analysis (see \cite{4}).
Our goal in this note is to illustrate this further by obtaining
Bourgain's result \cite{6} as a consequence of a generalized
version of M.\ Riesz's theorem due to 
Berkson, Gillespie, and Muhly \cite{5}.
 
We use the usual notation $L^p({\Omega},\mu,X)=\LBspc p$ to
denote the set of all strongly measurable functions on the
measure space $(\msp,\malg,\mu)$ with values in the Banach
space $X$ such that
$\myint_\Omega\|f(\omega)\|^p\,d\mu(\omega)<\infty$.  Recall
that $X$ is a UMD space if, for some (equivalently, all) $
p\in(1,\infty)$,
there exists a constant $C_p(X)$ such that for all $n\in
{\N}$,
$$ \big\| \sum_{j=1}^n \epsilon_j d_j \big\|_{L^p(X)}\leq
C_p(X)\big\| \sum_{j=1}^n d_j \big\|_{L^p(X)}\tag{1.1}$$
for every $X$-valued martingale difference sequence $(d_j)$
and for every $(\epsilon_j)\in\{-1,1\}^{\N}$. (For further 
background regarding this property, see \cite{7}, \cite{8} 
and \cite{9}.)
 
We say that $X$ is an HT space if for some
$p\in(1,\infty)$, there exists a constant $N_p(X)$ such that
for all $n \in \N$, we have
$ \|H_n f\|_{L^p(X)} \le N_p(X) \| f \|_{L^p(X)} $\ for all
$f \in L^p(\R,X)$, where
$$ H_n f(t) = {1\over \pi }\int_{1/n \le |s| \le n} {f(t-s)\over s}
\, ds \,.$$
 
Let $G$ be a compact abelian group with dual group $\Gamma$.
Let $\od\subset\Gamma$ be an order on $\dlgrp$, that is
$\od+\od\subset\od$, $\od\cap(-\od)=\{0 \}$, and
$\od\cup(-\od)=\Gamma$. Define a signum function $\sgn_\od$ on
$\dlgrp$ by
$\sgn_\od(\chi)=1,0$,or $-1$ according 
as $\chi\in\od\setminus\{0\}$,\ $\chi=0$, or  $\chi\in(-\od)\setminus\{0\}$.
Define a conjugate function operator $T_\od$ on the $X$-valued
trigonometric polynomials by
$$T_\od(\mysum_{\chi\in\dlgrp}a _\chi\,\chi) =-
i\mysum_{\chi\in\dlgrp}\sgn_\od(\chi) a _\chi\,\chi\,.
\tag{ 1.2}$$
Then we will say that a Banach space $X$ has the ACF
(abstract conjugate function) property if for some $p \in
(1,\infty)$, there is a constant $A_p(X)$ such that for all
compact abelian groups $G$\ with ordered dual groups, the operator
$T_\od$\ extends to $L^p(G,X)$, and  for 
all $f \in L^p(G,X)$, $ \|T_\od f\|_{L^p(X)}
\le A_p(X) \| f \|_{L^p(X)} $.
 
In \cite{8}, Burkholder showed that if $X$ is a UMD space,
it is an HT space.  He also conjectured that the converse was
true, and this was soon answered affirmatively by Bourgain
\cite{6}.  We will show the same result by first showing
that every HT space is an ACF space, and then showing that
every ACF space is a UMD space.
 
This new concept of ACF provides a natural bridge 
between the concepts of 
UMD and HT, thus demonstrating how the study of functions on
abstract abelian groups enhances and solidifies the connections
between harmonic analysis and martingale theory.  
 
\specialhead
2. Transference and the ACF Property
\endspecialhead
 
The following generalized version of M\. Riesz's theorem 
follows from \cite{5}, Theorem~4.1 
(see also \cite{2}, Theorem~6.3). 
We will sketch its proof to show the
role of transference and the HT property.  

\proclaim{(2.1) Theorem} If $X$ has the 
HT property, then it
has the ACF property.  Furthermore, for each $p \in
(1,\infty)$, we have that
$A_p(X) \le N_p(X)$.
\endproclaim
 
We follow the proof of 
 \cite{5}, Theorem~4.1.  We will 
need the following separation theorem for discrete 
groups (see \cite{3}, Theorem 5.14).  
As shown recently,
this result follows easily from the basic background 
required for the proof of Hahn's Embedding Theorem for 
orders (\cite{12}, Theorem 16, p.\ 16).  
For details, see \cite{4}.
 
\proclaim{(2.2) Theorem} Let $\od$ be an order on a discrete
abelian group $\dlgrp$. Given a finite subset
$F\subset\dlgrp$, there is a homomorphism $\psi:\dlgrp\to\R$
such that for all $\chi\in F$, $$\sgn_\od(\chi )=
\sgn(\psi(\chi ))\,.  \tag{2.2.1}$$
\endproclaim
 
\remark{(2.3) Proof of Theorem~2.1}
For $f$ an $X$-valued trigonometric polynomial, let $F$ be a 
finite subset of $\dlgrp$ such that
$f=\mysum_{\chi\in F} a_\chi\,\chi $ where $a_\chi\in X$.
Apply  Theorem~2.2 and obtain a real-valued homomorphism
$\psi $ such that $ \sgn(\psi(\chi ))=\sgn_\od(\chi ) $
for all $\chi\in F$. With this choice of $\psi$, $T_\od(f)=-i 
\mysum_{\chi\in F}\sgn(\psi(\chi))a_\chi\,\chi$.
 
Let $\phi:\R\to G$ be the adjoint homomorphism of $\psi$
which is defined by the relation
$\psi(\chi)(t)=\chi(\phi(t))$
for all $t\in\R$ and all $\chi\in\dlgrp$. It is easy to see
that for any
$\chi\in\dlgrp$, we have
$$\lim_{n\to\infty} {1\over \pi }\int_{1/n\le |t|\le n}{\chi(x-
\phi(t)) \over  t}\,dt=
     -i\,\sgn(\psi(\chi))\chi(x) \tag{2.3.1}$$
for all $x\in G$. Consequently, for an $X$-valued trigonometric 
polynomial, we have
$$\lim_{n\to\infty} {1\over \pi }\int_{1/n\le |t|\le n}{ f(x-
\phi(t)) \over  t}\,dt=
T_\od(f)(x)\tag{2.3.2}$$
for all $x\in G$. Thus, it is enough to
show that, for all $n\in\N$,
$$\bigg\| {1\over \pi}\int_{1/n\le |t|\le n}{f(x-\phi(t))\over
t}\,dt\bigg\|_{L^p(G,X)}\leq
N_p(X)\|f\|_{L^p(G,X)}\,.\tag{2.3.3}$$
 
This last inequality follows by adapting the transference argument 
of Calder\'on \cite{10}, 
Coifman and Weiss \cite{11}, 
to the setting of vector-valued functions,
and requires the HT property of $X$.
We omit its proof and 
refer the reader to the proof of Theorem~2.8 in \cite{5} 
for details.
\endprf\endremark
 
\specialhead
3. Proof of Bourgain's Result
\endspecialhead
 
As noted in \cite{9}, to show that $X\in\UMD$, it is enough
to consider dyadic martingale difference sequences defined
on $[0,1]$ using the Rademacher functions
$(r_n)_{n=1}^\infty$. (A proof of this reduction can be
obtained by adapting Remarque~3 in Maurey \cite{14} to the
vector-valued setting.) Also, as noted by Burkholder in  \cite{9}, 
it suffices to consider martingale difference sequences
such that $d_1=0$. In \cite{6}, Bourgain implicitly uses the
fact that to prove $X\in\UMD$, it is enough to consider
dyadic martingale difference sequences on the infinite
dimensional torus, $\T^\N$. We will give a precise version
of this, setting the notation in the process.
 
Let $\dualz=\wprod$ denote the weak direct product of $\N$ copies of $\Z$.
When endowed with the discrete topology, $\dualz$ is
topologically isomorphic to the dual group of $ \T^\N $. For
each $J=(j_n)\in\dualz$, we denote the corresponding
character by $\chi_J$, that is,
$\chi_J(\theta_1,\theta_2,\dots)=\prod\limits_{n=1}^\infty
e^{ij_n\theta_n}$ where, except for finitely many factors,
$e^{ij_n\theta_n}\equiv 1$. For each
$J=(j_n)\in\dualz\setminus\{ 0 \}$, define $n(J)$ to
be the largest $n\in\N$ such that $j_n\neq0$, and let
$n(0 )=0$.
With a slight abuse of notation,
we let $\Z^n$ denote $\{J\in\dualz:n(J)\le n \}$ for each
$n\ge 0$. Note that
$\dualz=\displaystyle{\cup_{n=1}^\infty\Z^n}$.
 
Identifying $\T$ with the interval $[-\pi,\pi)$, define a
sequence, $(s_n)_{n=1}^\infty$, of functions on $\T^\N$ by
$s_1 \equiv 1$, and for $n\ge 2$,
$s_n(\theta_1,\theta_2,\dots)=\sgn(\theta_{n-1})$.
Suppose  $(d_n)$ is a dyadic martingale difference sequence
on $[0,1]$ with $d_1=0$ and $d_n=v_n(r_1,\dots,r_{n-1})r_n$
where $v_n:\{-1,1\}^{n-1}\to X$ for all $n\ge 2$. Letting $d'_1=0$ and
$d_n^\prime=v_n(s_1,\dots,s_{n-1})s_n$ for $n\ge 2$ we
obtain a martingale difference sequence on $\T^\N$ such that
the sequences
$(d_n)$ and $(d_n^\prime)$ are identically distributed.
Call such a sequence a dyadic martingale difference
sequence on $\T^\N$. Therefore, to prove that $X\in\UMD$, it
suffices to show that there exists a constant satisfying
(1.1) for dyadic martingale difference sequences on $\T^\N$.
Since we may approximate each $d'_n$\ by a function with finite
spectrum in $\Z^n\setminus\Z^{n-1}$, to show that $X$\ has the 
UMD property, we see that it is
sufficient to show the following.
For $p \in (1,\infty)$,
there exists a constant $C_p(X)$ such that
if $K_j$ is a finite subset of $\Z^j\setminus\Z^{j-1}$, 
and $a_J\in X$ for all
$J\in K_j$, and
$(\epsilon_n)\in\{-
1,1\}^\N$, then
$$ \bigg\| \sum_{k=1}^n\epsilon_k\bigg(\sum_{J\in
K_j}a_J\,\chi_J
\bigg)\bigg\|_{L^p(\T^\N,X)}\leq
C_p (X)\bigg\| \sum_{j=1}^n \sum_{J\in K_j}
a_J\,\chi_J \bigg\|_{L^p(\T^\N,X)}\,.\tag{3.1}$$

This reduction appears in Bourgain \cite{6}.  We will show that if $X$ is 
an HT space, (3.1) follows from Theorem~2.1 with specific choices of the 
order $P$ on $\dualz$.  For this purpose, define a reversed lexicographic 
order on $\dualz$ as
follows:
$\od=\{J=(j_n)\in\dualz:j_{n(J)}>0\}\cup\{0 \}$
where as previously, if $J=(j_n)$, then $j_{n(J)}$ is the last
non-zero coordinate of
$J$. Thus, $\sgn_\od(\chi_J)=\sgn(j_{n(J)})$. Observe that
if $\epsilon=(\epsilon_n)\in\{-1,1\}^\N$,
then the set
$\od(\epsilon)=
\{J=(j_n)\in\dualz:\epsilon_{n(J)}j_{n(J)}>0\}\cup\{0 \}$
is also an order on $\dualz$. In this case, 
$\sgn_{\od(\epsilon)}(\chi_J)=\epsilon_{n(J)}\sgn(j_{n(J)})$.
We now state a simple identity that links the unconditionality of 
martingale difference sequences to harmonic conjugation with respect to 
orders:  for every $n\geq 1$, and all $J=(j_n)\in\Z^n\setminus\Z^{n-
1}$,\ we have
$$\epsilon_n  =\sgn_{\od(\epsilon)}(\chi_J)\,
\sgn_{\od}(\chi_J).\qquad
\tag{3.2}$$
To verify (3.2), simply note that for each $n\in\N$, 
if $J\in\Z^n\setminus\Z^{n-1}$, then $n(J)=n$.

From (3.2), one immediately obtains that
$$
T_P \circ T_{P(\epsilon)} \bigg( \sum_{k=1}^n \sum_{J\in K_k}
a_J\,\chi_J \bigg)
=
\sum_{k=1}^n\epsilon_k\bigg(\sum_{J\in
K_k}a_J\,\chi_J
\bigg),
\tag{3.3} $$
which expresses the martingale transform on the right side as a 
composition of two conjugate function operators. 
Applying Theorem~2.1 twice yields (3.1) and implies 
our next and last result.
 
\proclaim{(3.1) Theorem} Suppose $X$ is a Banach space, and
let $1<p<\infty$. If $X$ has ACF,
then $X$ is a UMD space and (1.1) holds with
$C_p(X) \leq (A_p(X) )^2$.\endproclaim
 
\remark{(3.2) Remarks} (a) Combining the implications above, we see that 
for a Banach space $X$, the properties
UMD, ACF, and HT are equivalent. 

(b) Burkholder also proved in \cite{8} that when $X$ is a UMD space,
the Hilbert transform is weak-type bounded on $L^1(\R,X)$.
This weak-type estimate also transfers to the conjugate function operator 
on $L^1(G,X)$ defined with respect to an arbitrary order on $\Gamma$.
For a proof see \cite{13}.
Thus, when $X\in\UMD$, we have vector-valued analogs of
the classical results of M\. Riesz and Kolmogorov. 

One may ask for an analog of Privalov's classical theorem for 
UMD-valued functions in $L^1(G,X)$. For the case of $X=\C$, this is done 
in \cite{4}. 
 
Transference of operators on vector-valued
$L^p$-spaces can be carried out in much greater generality
than that shown here. A representation satisfying the
vector-valued version of the distributional control
condition introduced in \cite{1} will transfer strong-type
and weak-type bounds for maximal operators. This work will
appear in \cite{13}. The contents of this article are part
of the second author's dissertation.
\endremark

Department of Mathematics, University of Missouri-Columbia,
Columbia, Missouri
 
\noindent 65211 U\.S\.A.
 
\enddocument